
\documentstyle{amsppt}
\pagewidth{6 in}\vsize8in\parindent=6mm\parskip=3pt\baselineskip=14pt
\tolerance=10000\hbadness=500\loadbold
\magnification=1100
\topmatter
\title
Oscillatory  integral operators with low--order degeneracies
\endtitle
\author
Allan Greenleaf and Andreas Seeger
\endauthor
\abstract
We prove sharp $L^2$ estimates for oscillatory integral and Fourier integral 
operators  for which the 
associated canonical relation $\Cal C\subset T^*\Omega_L\times T^*\Omega_R$
projects to $T^*\Omega_L$ and to $T^*\Omega_R$ with corank one singularities 
of type $\le 2$. This includes two-sided cusp singularities.
Applications are given to operators with one-sided swallowtail singularities 
such as restricted X-ray transforms for well-curved line complexes in five dimensions.
\endabstract
\thanks
Research supported in part by NSF grants DMS 9877101 (A.G.) and 
DMS 9970042 (A.S.).
\endthanks
\address
University of Rochester, Rochester, NY 14627
\endaddress
\address
University of Wisconsin, Madison, WI 53706
\endaddress
\leftheadtext{Allan Greenleaf and Andreas Seeger}
\rightheadtext{Integral operators with low-order degeneracies}
\endtopmatter
\document

\define\ic{\imath}
\define\vth{\vartheta}
\define\ffx{x}
\define\ffz{z}

\define\cf{{\it cf}}

\define\rank{{\text{\rm rank }}}

\define\loc{{\text{\rm loc}}}

\define\comp{{\text{\rm comp}}}

\define\supp{{\text{\rm supp }}}

\define\inn#1#2{\langle#1,#2\rangle}

\define\noi{\noindent}
\define\lcontr{\rfloor}
\define\lco#1#2{{#1}\lcontr{#2}}
\define\lcoi#1#2{\imath({#1}){#2}}
\define\rco#1#2{{#1}\rcontr{#2}}
\redefine\exp{{\text{\rm exp}}}
\define\bin#1#2{{\pmatrix {#1}\\{#2}\endpmatrix}}

\define\lc{\lesssim}
\define\gc{\gtrsim}


\define\eps{\varepsilon}

\define\ka{\kappa}
            
\define\la{\lambda}             
             
\define\si{\sigma}

\define\om{\omega}              \define\Om{\Omega}


\define\bbR{{\Bbb R}}

\define\bbZ{{\Bbb Z}}

\define\cA{{\Cal A}}

\define\cC{{\Cal C}}

\define\cH{{\Cal H}}

\define\cK{{\Cal K}}
\define\cL{{\Cal L}}

\define\cR{{\Cal R}}

\define\cT{{\Cal T}}
\define\cU{{\Cal U}}
\define\cV{{\Cal V}}

\define\cZ{{\Cal Z}}





\define\set#1{\{{#1}\}}
\define\emph#1{{{\it #1}}}
\define\norm#1{\|{#1}\|}
\define\bignorm#1{\big\|{#1}\big\|}
\define\Bignorm#1{\Big\|{#1}\Big\|}
\define\spanset#1{{\text{span}\{#1\}}}

\head{\bf Introduction}\endhead
Let $\Omega_L$, $\Omega_R$ be open sets in $\Bbb  R^d$.
This paper is concerned with $L^2$ bounds for
 oscillatory integral operators $T_\la$ of the form
$$T_\la f(x)=\int e^{\ic \la\Phi(x,z)} \si(x,z) f(z) dz
\tag 1.1$$
where $\Phi\in C^\infty(\Om_L\times\Om_R)$ is real-valued,
 $\si\in C^\infty_0(\Omega_L\times\Omega_R)$ and $\la$ is large.
We shall also write
$$T_\la\equiv T_\la[\si]
$$
to indicate the dependence on the symbol $\si$.

The decay in $\la$ of the $L^2$ operator norm of $T_{\la}$ is  determined by
the geometry of the canonical relation
$$\cC=\{(x,\Phi_x,z,-\Phi_z): (x,z)\in \Omega_L\times\Omega_R\}\subset
T^*\Omega_L\times T^*\Omega_R,
\tag 1.2$$
specifically by the behavior of the projections
$\pi_L:\cC\to T^*\Om_L$ and $\pi_R:\cC\to T^*\Om_R$ ,
$$ \aligned
&\pi_L:(x,z)\mapsto (x,\Phi_x(x,z))
\\
&\pi_R:(x,z)\mapsto (z,-\Phi_z(x,z));
\endaligned
\tag 1.3
$$
here $\Phi_x$ and $\Phi_z$ denote the partial gradients with respect to $x$
and $z$.
Note that  $\rank D\pi_L=\rank D\pi_R$ is equal to $d+\rank\Phi_{xz}$ and
that
the determinants of $D\pi_L$ and $D\pi_R$ are equal to
$$h(x,z):= \det \Phi_{xz}(x,z).
\tag  1.4
$$

If $\cC$ is locally the graph of a canonical transformation, i.e., if 
$h\neq 0$,  then
$\|T_\la\|=O(\la^{-d/2})$ (see H\"ormander \cite{15}, \cite{16}). If the
projections have singularities
then there is less decay in $\la$ and  in various specific cases
the decay has been determined. In dimension  $d=1$
Phong and Stein \cite{21} obtained a
complete description of the $L^2$ mapping properties, for the case of
real-analytic
phase functions.
Similar results for $C^\infty $ phases  (which however  missed  the
endpoints)
 and related $L^p$ estimates for  averaging operators in the plane are in
\cite{24}. The bounds  for oscillatory integral operators  in one dimension,
with
$C^\infty$ phases, have recently
been substantially improved by Rychkov \cite{22}, so that
many endpoint estimates  are now available in the $C^\infty$ category.

Such general results are not known in higher dimensions even under the
assumption of $\rank
\Phi_{xz}\ge  d-1$. We list some known cases.
If both projections $\pi_L$ and $\pi_R$ have fold $(S_{1,0})$ singularities
then
$\|T_\la\|=O(\la^{-(d-1)/2-1/3})$ (\cite{17}, \cite{19}, \cite{5}).
If only one  of the projections  has fold singularities
then by \cite{8} we have
$\|T_\la\|=O(\la^{-(d-1)/2-1/4})$; this is sharp if the other projection
is maximally degenerate (\cite{13}) but can be improved when that projection
satisfies some finite type finite type condition (for sharp results of this
sort see Comech \cite{3}).
This one-sided behavior comes up naturally when studying  restricted X-ray
transforms \cite{6}, \cite{11}, \cite{14}. In \cite{9}
 the authors began a study of the case
of higher one-sided Morin $(S_{1_r,0})$ singularities, which are the stable
singularities of corank one,
and
 it was shown under suitable additional
(''strongness'') assumptions  that such estimates can be deduced
from sharp estimates for two-sided $S_{1_{r-1},0}$ singularities. Thus the
authors were able to prove that if one projection is a Whitney cusp,
i.e., of type $S_{1,1,0}$, then  $\|T_\la\|=O(\la^{-(d-1)/2-1/6})$;
again  this is only
 sharp if the other projection is maximally degenerate.

It is conjectured that if one of $\pi_L$ or $\pi_R$ has $S_{1_r,0}$
singularities then
$\|T_\la\|=O(\la^{- (d-1)/2-1/(2r+2)})$ (for the discussion
of some model cases where this is satisfied and sharp see \cite{9}).
Here we take up the case $r=3$; such mappings are commonly
referred to as  {\it swallowtail} singularities. In order to prove this
result  it is  crucial to get a sharp  result for
operators with  two-sided cusp singularities.

\proclaim{Theorem}

(i)  Suppose that the only singularities of one of the projections
($\pi_L$ or $\pi_R$)  are Whitney folds, Whitney cusps or swallowtails.
Then  $\|T_\la\|_{L^2\to L^2}=O( \la^{-(d-1)/2 -1/8})$ for $\la\ge 1$.

(ii)
Suppose that the only singularities of both  projections
$\pi_L$ and $\pi_R$  are Whitney folds or  Whitney cusps. Then
$\|T_\la\|_{L^2\to L^2}=O( \la^{-(d-1)/2-1/4})$ for $\la\ge 1$.
\endproclaim

A slightly weaker result than (ii)  was recently obtained by
Comech and Cuccagna \cite{4}, who proved for two-sided cusp singularities the
bound
$\|T_\la\|_{L^2\to L^2}\le C_\eps \la^{-(d-1)/2 -1/4+\eps}$
with $C_\eps\to \infty$
 as $\eps\to 0$.

We shall prove somewhat
more general results about operators of the same ``type'' but with
the stability assumptions weakened.
To formulate the hypotheses  we review the definition of kernel vector
fields
for a map.
Fix $n$-dimensional manifolds $M , N $ and points $P_0\in M $ and
$Q_0\in N $.
Let $f:M \to N $ be a $C^\infty$ map with $f(P_0)=Q_0$.
Let $\cU$ be a neighborhood of $P$.
A vector field $V$ is a {\it kernel  field} for the map $f$ on $\cU$  if $V$
is smooth on $\cU$
 and
if   $Df_P V= \det (Df_P)W_{f(P)}$ for $P\in \cU$; here  $W$ is a  smooth
vector field on $N $  defined near  $Q_0=f(P_0)$ and $\det (Df_P)$ is
calculated with respect to any local systems of coordinates.

Suppose now that  $\rank Df(P_0)\ge n-1$. Then there  is a neighborhood of
$P$ and a nonvanishing kernel vector field $V$ for $f$ on $\cU$. If
$\widetilde V$ is another kernel field on $U$
 then $\widetilde V= \alpha V-\det(Df)W$
in some neighborhood of $P_0$, for some vector field $W$ and smooth function
$\alpha$.

This is easy to see by an elementary calculation. Indeed we
may choose coordinates $x=(x',x_n)$ on $M $, $y=(y',y_n)$ on $N $
 vanishing at $P_0$ and $Q_0$, respectively, so that
$D_x'f=(A,b)$ and $D_{x_n} f= (c^t,d) $
where $A$ is an invertible $(n-1)\times (n-1)$ matrix, $b$ and $c$
 are vectors in $d\in \bbR$, and  $A,b,c,d$ depend smoothly on $x$.
Define the vector field $V$ by
$V=\partial_{ x_n}-\inn {A^{-1} b}{\partial_{x'}}$.
Then clearly $Df(V)=(d-c^tA^{-1}b)\partial_{y_n}$
and
$\det Df= (d-c^tA^{-1}b) \det A$; thus $V$ 
is a kernel field. Now assume that
 $\widetilde V= \inn{\beta'}{\partial_{x'}}+\beta_n
\partial_{x_n}
$
so that  $Df(\widetilde V)=\det (Df )Z$ with
$Z=Z'+\gamma_n\partial_{y_n}$, $Z'=\inn{\gamma'}{\partial_{y'}}$;
here $\beta=(\beta',\beta_n)$ are smooth functions of $x$ and
 $\gamma=(\gamma',\gamma_n)$ are smooth functions of $y$. Then,
at any $x$,
 $A\beta'+b\beta_n=\det (Df) \beta'$; therefore
$\beta'=-A^{-1}b \beta_n+\det Df A^{-1}\gamma'$ and  thus
$
\widetilde V=\beta_n V+(\det Df) Z'
$ as claimed.

\definition{Definition}  Suppose that
$M $ and $N $ are smooth $n$-dimensional manifolds and
that
 $f:M \to N $ is a smooth map with   $\dim \text{ker}( Df) \le 1$
on $M $.
We say that $f$ is of type $k$ at $P$ if
there is a nonvanishing kernel field $V$ near $P$
 so that $V^j(\det Df)_P=0$ for $j<k$ but
$V^k(\det Df)_P\neq 0$.
\enddefinition

>From the previous discussion it is clear that this  definition does not
depend on the choice  of the nonvanishing kernel field. If one assumes
that $Df$ drops rank {\it simply} on the singular variety $\{ \det
Df=0\}$  (i.e.,  if $\nabla\det Df\neq 0$) then the  definition
agrees with the one proposed by Comech  \cite{3}.

\proclaim{Theorem 1.1} Suppose that both $\pi_L$ and $\pi_R$ are of type
$\le 2$
on $\cC$. Then  for $\la\ge 1$
$$\big\|T_\la\big\|_{L^2\to L^2}=O(\la^{-(d-1)/2-1/4}).
$$
\endproclaim

\proclaim{Theorem 1.2} Suppose that  $D\pi_L$ drops rank simply on the
singular variety
$\{\det D\pi_L=0\}$  and suppose that  $\pi_L$ is  of type $\le 3$
on $\cC$. Then  for $\la\ge 1$
$$\big\|T_\la\big\|_{L^2\to L^2} =O(\la^{-(d-1)/2-1/8}).
$$
\endproclaim

Of course the analogous statement holds with $\pi_L$ replaced by $\pi_R$ in
Theorem 1.2.
As a corollary of both theorems we obtain the sharp endpoint
estimate for two-sided cusp and one-sided swallowtail singularities stated
above.

\remark{\bf Remark}
The estimates in Theorem 1.1 and Theorem 1.2 are
 stable under small perturbations of $\Phi$ and $\si$ in the
$C^{\infty}$-topology.
\endremark

 The above theorems imply sharp $L^2$-Sobolev estimates
for Fourier integral operators (see \cite{8}).
Let
$C\subset T^*\Om_L\setminus\{0\}\times T^*\Om_R\setminus\{0\}$
and let  $F\in I^\mu(\Omega_L,\Omega_R; C)$ (see \cite{15} for the
definition
and \cite{8} for the reduction of smoothing estimates for Fourier
integral operators to decay estimates for oscillatory integral operators).
As a corollary of Theorems 1.1 and 1.2 one obtains

\proclaim{Theorem 1.3}

(i) If both $\pi_L$ and $\pi_R$
are  of type $\le 2$, then $F$ maps  $L^2_{\alpha,\comp}$
to $L^2_{\alpha-\mu-1/4,\loc}$.

(ii) If one projection ($\pi_L$ or $\pi_R$) is of type $\le 3$ and the
rank of its differential  drops only simply, then  $F$ maps
$L^2_{\alpha,\comp}$
to $ L^2_{\alpha-\mu-3/8,\loc}$. \endproclaim

\medskip
\remark{\bf Remarks} {\bf   1.} As an example of part (i) of Theorem
1.3, consider as in \cite{9,\S6} a family of curves in $\Bbb R^4$ of the
form $$\gamma_x(t)=\exp(tX+t^2Y+t^3Z+t^4W)(x)$$
for smooth vector fields $X,Y,Z,W$ on $\Bbb R^4$ such that both of the
sets of vectors
$$\Bigl\{X,Y,Z\mp\frac16[X,Y],W\mp\frac14[X,Z]+\frac1{24}[X,[X,Y]]\Bigr\}$$
are
linearly independent at each point $x$. Then the generalized Radon
transform $$Rf(x)=\int_{\Bbb R} f(\gamma_x(t))\chi(t) dt,\quad\chi\in
C_0^\infty(\Bbb R),$$ belongs to $I^{-\frac12}(\Bbb R^4,\Bbb R^4;C)$ with
the canonical relation $C$ a two--sided cusp, i.e., both $\pi_L$ and
$\pi_R$ are Whitney cusps, and hence it follows from Theorem 1.3(i) that
$R:L^2_{\alpha,\comp}\rightarrow L^2_{\alpha+1/4,\loc}$,
for all $\alpha\in\Bbb R$, generalizing the well-known fact for the
translation--invariant family $\gamma_x(t)=x+(t,t^2,t^3,t^4)$.

{\bf 2.} Consider the translation--invariant families of curves in $\Bbb
R^3$, $\gamma_x^1(t)=x+(t,t^2,t^4)$ and $\gamma_x^2(t)=x+(t,t^3,t^4)$.
Then $\{\gamma_x^1\}$ is associated with a canonical relation, $C^1$, which
is a two--sided cusp, while $\{\gamma_x^2\}$ is associated with a canonical
relation, $C^2$, for which both projections are type 2, but not Whitney
cusps.
In fact, the singular variety of $C^2$ is not smooth: it is  a union
of two intersecting hypersurfaces, and   $det(D\pi_L)$ and $det(D\pi_R)$
vanish
of order two at  the intersection and simply away from that
intersection. Averaging operators associated with
any (not necessarily translation--invariant) sufficiently small $C^\infty$
perturbation of either $\{\gamma_x^1\}$ or  $\{\gamma_x^2\}$
 will still have both projections of type 2 and
hence map $L^2_{\alpha,\comp}\rightarrow L^2_{\alpha+\frac14,\loc}$.

{\bf 3.} As an instance of part  (ii) of Theorem 1.3,  let $\cR$ be
 the restricted X-ray transform associated to a
well--curved line complex $\frak C$ in $\bbR^5$
(see \cite{9,\S5} for the  definition). Then $\pi_R$ has (at most)
swallowtail singularities and  $\cR$ maps
$L^2_\comp(\bbR^5)$ into  $L^2_{1/8,\loc}(\frak C)$.
As an example consider a curve $\alpha\to \gamma(\alpha)$ in
$\Bbb R^4$
with
$\gamma'$, $\gamma''$, $\gamma'''$ and $\gamma^{(4)}$  being linearly
independent at each $\alpha$ and consider
the X-ray transform associated to the rigid $5$-dimensional line complex
consisting of lines
$\{\ell_{x',\alpha}: x'\in \bbR^4, \alpha\in \bbR\}$ in $\bbR^5$ where
$\ell_{x',\alpha}=\{(x'+t\gamma(\alpha),t), t\in \bbR\}$, and perturbations
of
this example. For the rigid  case the projection $\pi_L$ is a blowdown
in the sense of \cite{13} or \cite{14},
i.e., it exhibits a maximal degeneracy; this behavior however  is not
invariant under small perturbations and is not required for Theorem 1.3 to
apply.

{\bf 4.} As an example of a restricted $X$-ray transform in $\bbR^4$
which is not well--curved in the sense of \cite{9}, consider the situation
as in the previous example, but with $\gamma$ replaced by one of the
curves   $\gamma^{(1)}(\alpha)=(\alpha,\alpha^2, \alpha^4)$ or
$\gamma^{(2)}(\alpha)= (\alpha,\alpha^3, \alpha^4)$ in $\bbR^3$. For both
examples $\pi_R$  satisfies a type three condition with $\det d\pi_R$
vanishing  simply; however  the singularity of  $d\pi_R$ for the
canonical relation associated to the second line complex (defined by
$\gamma^{(2)}$) is not of swallowtail type. Again $\cR$ and perturbations
thereof map $L^2_\comp(\bbR^4)$ into  $L^2_{1/8,\loc}(\frak C)$.

{\bf 5.}
For conormal operators in two dimension the   condition  of type $\le k$
for $\pi_L$  corresponds to a left
finite type condition of order $k+2$ in the terminology of \cite{23}, and
the condition of (exact) type $k$ corresponds to the type  $(1,k+1)$
condition
in the terminology of \cite{24}.
\endremark

\head{\bf 2. Bounds for operators with two-sided type two conditions}
\endhead

We decompose
the operator according to the size of
$\det \Phi_{xz}$, following  Phong and Stein \cite{20} who used this
decomposition  to estimate
 operators with fold singularities.
 Various extensions and refinements
are in \cite{23}, \cite{21}, \cite{5},  \cite{10}, \cite{3}, \cite{4};
in fact we will use the
 key estimate in \cite{4} as the first step in our proof of Theorem 1.1.
As in  that work (see also \cite{23}, \cite{3})  we shall need to localize
$V_Lh$ and $V_Rh$
where $V_L$ and $V_R$ are nonvanishing kernel vector fields for $\pi_L$ and
$\pi_R$, respectively.
We may suppose that the support of $\si$ is small and choose
 coordinates $x=(x',x_d)$, $z=(z',z_d)$ in $\bbR^{d-1}\times\bbR$
vanishing at a reference point $P^0=(x^0,z^0)$ so that
$$\Phi_{x'z'}(P^0)=I_{d-1}, \qquad
\Phi_{x_dz'}(P^0)=0,\qquad
\Phi_{x'z_d}(P^0)=0.
$$
Write $\Phi^{z'x'}:=\Phi_{x'z'}^{-1}$ and
$\Phi^{x'z'}:=\Phi_{z'x'}^{-1}= (\Phi_{x'z'}^t)^{-1}$.
Representatives for the kernel vector fields are then given by
$$\aligned
V_R&=\partial_{x_d}-\Phi_{x_dz'}\Phi^{z'x'}\partial_{x'}
\\
V_L&=\partial_{z_d}-\Phi_{z_dx'}\Phi^{x'z'}\partial_{z'}
\endaligned
\tag 2.1
$$
(see \cite{2} and the discussion in the introduction).

Let $K$ be a fixed compact set in $\Omega_L\times\Omega_R$ which contains
the support of $\sigma$ in its interior.
Let $A_0\ge 10^{2d}$ so that
$$
\|\Phi\|_{C^5(K)}\le 10^{-2d} A_0.
\tag 2.2
$$
We also assume that
$$
 |V_L^2h|\ge A_1^{-1}, \ \
|V_R^2 h|\ge A_1^{-1}.
\tag 2.3
$$
for some  $A_1\ge 1$.
After additional localization we may assume that $\sigma$ is supported on a
set of
small diameter $\eps$, for later use we choose
$$\eps= 10^{-1}\min\{A_0^{-2}, A_1^{-2}\}.
\tag 2.4$$

Let $\beta_0\in C^\infty(\Bbb R)$ be an even function supported in $(-1,1)$,
and equal to one in $(-1/2,1/2)$. Let $\beta(s)\equiv \beta_1(s)=
\beta_0(s/2)-\beta_0(s)$ and for
$j\ge 1$ let
$\beta_j(s)=\beta_1(2^{1-j}s)=
\beta_0(2^{-j}s)-\beta_0(2^{-j+1}s)$.

We may assume that $\la$ is large.
Let $\ell_0=[\log_2(\sqrt\la)]$, that is  the largest integer $\ell$
 so that $2^{\ell}\le \la^{1/2}$.
Let  then
$$
\aligned
\si_{j,k,l}(x,z)&= \si(x,z) \beta(2^l h(x,z))\beta_j( 2^{l/2}V_Rh(x,z))
\beta_k(2^{l/2}V_Lh(x,z))
\\
\si^0_{j,k,\ell_0}(x,z)&=
 \si(x,z) \beta_0(2^\ell_0 h(x,z))\beta_j( 2^{\ell_0/2}V_Rh(x,z))
\beta_k(2^{\ell_0/2}V_Lh(x,z));
\endaligned
\tag 2.5
$$
thus if $j,k>0$ then
 on the support of $\si_{j,k,l}$ we have that $|h|\approx 2^{-l}$,
$|V_L h|\approx 2^{k-l/2}$,
$|V_R h|\approx 2^{j-l/2}$.

Our main technical result sharpens estimates given in \cite{4}; we
use here, as throughout, the notation $A\lc B$ to denote
inequalities $A\le CB$ with constants  $C$ independent of
$\lambda, j,k, l$.

\proclaim{Theorem 2.1}
We have the following estimates:

(i) For $0<l<\ell_0=[\log_2(\sqrt\la)]$
$$
\|T_\la[\si_{j,k,l}]\|_{L^2\to L^2}\lc
\la^{-(d-1)/2}\min\big \{ 2^{l/2}
\la^{-1/2}; 2^{-(l+j+k)/2}
\big\}.
\tag 2.6
$$

(ii)
$$
\|T_\la[\si^0_{j,k,\ell_0}]\|_{L^2\to L^2}\lc
\la^{-(d-1)/2-1/4} 2^{-(j+k)/2}.
\tag 2.7
$$
\endproclaim

Given Theorem 2.1 we can deduce Theorem 1.2 by
 simply summing  the estimates (2.6) and (2.7):
The bound
$\sum_{j,k}\|T_\la[\sigma_{j,k,\ell_0}^0]\|\le \la^{-(d-1)/2-1/4}$ is immediate.
Moreover
$$
\sum_{0\le l\le \log_2(\sqrt\la)}\sum_{0\le j,k\le l/2}
\|T_\la[\sigma_{j,k,l}]\|_2\le I+II
$$
where
$$
\align
I&\le
\sum_{0\le l\le \log_2(\sqrt\la)}
\sum\Sb j,k \\
j+k\le  \log_2(\la 2^{-2l})\endSb 2^{l/2}\la^{-d/2}
\\
&\lc \la^{-(d-1)/2-1/4}
\sum_{0\le l\le \log_2(\sqrt\la)} (\la 2^{-2l})^{-1/4} [1+\log(\la 2^{-2l})]^2
\lc  \la^{-(d-1)/2-1/4}
\endalign
$$
and
$$
\align
II&\le
\sum_{0\le l\le \log_2(\sqrt\la)} 2^{-l/2}\la^{-(d-1)/2}
\sum\Sb j,k \\
j+k\ge  \log_2(\la 2^{-2l})\endSb
 2^{-(j+k)/2}
\\
&\lc \la^{-(d-1)/2}
\sum_{0\le l\le \log_2(\sqrt\la)} 2^{-l/2} (2^l\la^{-1/2})(1+ \log_2(\la
2^{-2l}))
\lc  \la^{-(d-1)/2-1/4}. \qed
\endalign
$$

We make some preliminary observations needed in the proof of Theorem 2.1.
In what follows we always   make the

\noi{\bf Assumption:}  $k\le j$.

\noi For  $k\ge j$ apply the corresponding  estimates
for the adjoint of $T_\la[\si_{j,k,l}]$.
For the proof of Theorem 2.1 we may assume, by the known result for
one-sided folds
\cite{8}, that
$$ 2^{k-l/2}\le 2^{j-l/2}\le \eps
\tag 2.8
$$
where $\eps$ is as in (2.4).

\subhead{Affine changes of variables}\endsubhead
Before starting with estimates we wish to mention the effect of changes of
variables  on  (2.5). Set $x=\ffx(u)$ and $z=\ffz(v)$ and let
$\Psi(u,v)=\Phi(\ffx(u), \ffz(v))$. Let
$h(x,z)=\det \Phi_{xz}$ and $\widetilde h(u,v)= \det\Psi_{uv}(u,v)$ then
$$\widetilde  h(u,v) = h(\ffx(u), \ffz(v))
\det \frac{D\ffx}{Du}
\det \frac{D\ffz}{Dv}.
$$
If
$V_R=\sum s_i(x,z)\partial_{x_i}$,
$V_L=\sum t_i(x,z)\partial_{z_i}$
and
$\widetilde V_R=\sum \sigma_i(u,v)\partial_{u_i}$,
$\widetilde V_L=\sum  \tau_i(u,v)\partial_{v_i}$,
then
$V_R g(x(u),z(v))= \widetilde V_R(g(x(u),z(v))$ and
$V_L g(x(u),z(v))= \widetilde V_L(g(x(u),z(v))$ if and only if
 $\vec \sigma = (\tfrac{D\ffx}{Du})^{-1} \vec s$ and
$\vec \tau = (\tfrac{D\ffz}{Dv})^{-1} \vec t$.

In particular if our changes of variables are affine   and of the form
$$
\ffx(u)=x^0+(u'+ a' u_d,u_d) , \
\ffz(v)=z^0+(v'+ b' v_d,v_d)
\tag 2.9
$$ with constant vectors $a',b'\in\bbR^{d-1}$ and $P=(x^0,z^0)$ and if
 $V_L$ and $V_R$ are of the form (2.1) then we have
$\widetilde V_R=\partial_{u_d}-(a^t+\Phi_{x_dz'}\Phi^{z'x'})\partial_{u'}$
and
$\widetilde V_L=\partial_{v_d}-(b^t+\Phi_{z_dx'}\Phi^{x'z'})\partial_{v'}$
where all coefficient  functions are evaluated at
$(x^0,z^0)+(u'+ a' u_d,u_d,v'+ b' v_d,v_d)$.
Thus by choosing
$ a'=- \Phi^{x'z'}(P) \Phi_{z'x_d}(P)$, $b'=-\Phi^{z'x'}(P)\Phi_{x'z_d}(P)$
we  achieve that 
$(\widetilde V_R)_{0,0}=\partial_{u_d}$,
$(\widetilde V_L)_{0,0}=\partial_{v_d}$.

\subhead{Localization}\endsubhead
We shall perform various localizations to small boxes in $(x,z)$-space.

Let $P=(x^0,z^0)\in \Omega_L\times \Omega_R\subset \Bbb R^d_L\times \Bbb
R^d_R$ and let
$a\in \Bbb R^{d}_L$ and  $b\in \bbR^d_R$ be  vectors with $1\le
\|a\|_\infty,
\|b\|_\infty\le 2$
 and let $\pi_a^\perp$, $\pi_b^\perp$ be the orthogonal projections  to
the orthogonal complement of $\bbR a$ in $\bbR^d_L$ and
$\bbR b$ in $\bbR^d_R$,
respectively. Suppose  $0<\gamma_1\le \gamma_2<1$ and
$0<\delta_1\le\delta_2\le \eps$   and
let
$$\multline
B_P^{a,b}(\gamma_1,\gamma_2,\delta_1,\delta_2)
=\\
\{(x,z): |\pi_a^\perp(x-x^0)|\le \gamma_1, |\inn{x-x^0}{a}|\le \gamma_2,
|\pi_b^\perp(z-z^0)|\le \delta_1, |\inn{z-z^0}{b}|\le \delta_2.
\}\endmultline
\tag 2.10
$$

\definition{Definition}
We say that $\chi\in C^\infty_0$ is a {\it  normalized cutoff function
associated to $B_P^{a,b}(\gamma_1,\gamma_2,\delta_1,\delta_2)$} if it is
supported in
$B_P^{a,b}(\gamma_1,\gamma_2,\delta_1,\delta_2)$ and satisfies the (natural)
estimates
$$
|(\pi_a^\perp\nabla_x)^{m_L}
\inn{a}{\nabla_x}^{n_L}
(\pi_b^\perp\nabla_z)^{m_R}
\inn{b}{\nabla_z}^{n_R}
\chi(x,z)|\le  \gamma_1^{-m_L}\gamma_2^{-n_L}
 \delta_1^{-m_R}\delta_2^{-n_R}
$$
whenever  $m_L+n_L\le 10d$,  $m_R+n_L\le 10d$.
Here
$(\pi_a^\perp\nabla_x)^{n_L}$ stands for any differential operator
$\inn{\vec u_1}{\partial_x}\dots
\inn{\vec u_{n_L}}{\partial_x}
$
where the vectors $u_1,\dots, u_{n_L}$ are unit vectors perpendicular to
$a$.

We denote by $\cZ_P^{a,b}(\gamma_1,\gamma_2,\delta_1,\delta_2)$ the class of
all normalized cutoff functions associated to
$B_P^{a,b}(\gamma_1,\gamma_2,\delta_1,\delta_2)$.

\enddefinition

We shall often localize to boxes of the form (2.10) and consider
$T_\la[\zeta \si_{j,k,l}]$ where $\zeta$ is a
cutoff function which is  controlled
by an absolute constant times a normalized cutoff function in the above
sense.

Suppose now that $P=(x^0,z^0)$
and  our change of variable is as in (2.9) and that $a=(a',1)$,
$b=(b',1)$. Suppose that $\zeta$ is a normalized cutoff
function associated  to
$B_P^{a,b}(\gamma_1,\gamma_2,\delta_1,\delta_2)$.
Let
 $\widetilde \zeta(u,v)=\zeta(\ffx(u),\ffz(v))$. Then
$\widetilde \zeta$  is supported in
$\widetilde B_{(0,0)}^{e_d, e_d}(\tilde\gamma_1,\tilde\gamma_2,
\tilde\delta_1,\tilde
\delta_2)$ with
$\tilde\gamma=(1+|a'|)\gamma$,
$\tilde \delta=(1+|b'|)\delta$ and
there is a  positive  constant $C$ (independent of $\gamma$, $\delta$) so
that
$C^{-1}\widetilde \zeta$
is a normalized cutoff function associated to
$\widetilde B_{(0,0)}^{e_d, e_d}(\tilde\gamma_1,\tilde\gamma_2,
\tilde\delta_1,\tilde
\delta_2)$.

Changing variables as in (2.9)
in the expression for the  operator $T_\la[\zeta \si_{j,k,l}]$ yields that
$$T_\la[\zeta \si_{j,k,l}]f(\ffz(v))= \int \widetilde \zeta(u,v) \widetilde
\si_{j,k,l}
e^{\ic \la\Psi(u,v)} dv
$$
with
$
\widetilde \si_{j,k,l}(u,v)= \si(\ffx(u),\ffz(v)) \beta(2^l\widetilde
h(u,v))
\beta_j( 2^{l/2}\widetilde V_R\widetilde h(u,v))
\beta_k(2^{l/2}\widetilde V_L\widetilde h(u,v))$
and $\widetilde h(u,v)=\det\Psi_{uv}(u,v)$.

\subhead{Basic estimates}\endsubhead

We now give estimates for various pieces localized to (thin) boxes which
will usually
be  longer in the directions of the kernel  fields $V_R$ and $V_L$.

In order to formulate our results we start with a definition.

\definition{Definition}
 Let $P=(x^0,z^0)\in\Omega_L\times \Omega_R$ and let
$$ a_P=\big(
- \Phi^{x'z'}(P) \Phi_{z'x_d}(P),1\big), \qquad
b_P=\big(-\Phi^{z'x'}(P)\Phi_{x'z_d}(P),1\big)
\tag 2.11
$$
Define, for fixed $j,k,l$,
$$
\align
&\cA_P(\gamma_1,\gamma_2,\delta_1,\delta_2):=\sup\big\{
\big\|T_\la[\zeta \si_{j,k,l}]\big\|: \ \zeta\in
\cZ_P^{a_P,b_P}(\gamma_1,\gamma_2,\delta_1,\delta_2)\big\}
\tag 2.12
\\
&\cA_P^0(\gamma_1,\gamma_2,\delta_1,\delta_2):=\sup\big\{
\big\|T_\la[\zeta \si^0_{j,k,\ell_0}]\big\|: \ \zeta\in
\cZ_P^{a_P,b_P}(\gamma_1,\gamma_2,\delta_1,\delta_2)\big\}.
\tag 2.13
\endalign$$
Here $\ell_0=[\log_2(\sqrt\la)]$.
\enddefinition

The main estimate   in  Comech-Cuccagna\cite{4} applies to
 operators whose kernels
are localized to boxes
$B_P^{a_P,b_P}
(2^{-l}, 2^{-j-l/2}, 2^{-l}, 2^{-k-l/2})$. This
result is formulated in (2.14) of the
following proposition. The constants implicit in the inequalities below,
do not depend on $j,k,l$.

\proclaim{Proposition 2.2 }
(i) For $2^l\le \la^{1/2}$, $k\le j\le l/2$,
$$
\sup_P \, \cA_P(2^{-l}, 2^{-j-l/2}, 2^{-l}, 2^{-k-l/2})
\lc 2^{l/2}\la^{-d/2}.
\tag 2.14
$$
and
$$
\sup_P \, \cA_P(2^{-l}, 2^{-j-l/2}, 2^{-l}, 2^{-k-l/2})
\lc 2^{-(l+j+k)/2}\la^{-(d-1)/2}.
\tag 2.15
$$

(ii)  Let
$l=[\log_2(\sqrt\la)]=\ell_0$. Then for
$k\le j\le \ell_0/2$,
$$
\sup_P \, \cA_P^0(2^{-\ell_0}, 2^{-j-\ell_0/2},
 2^{-\ell_0}, 2^{-k-\ell_0/2})
\lc
\la^{-(d-1)/2-1/4} 2^{-(j+k)/2}.
\tag 2.16
$$
\endproclaim

Proposition 2.2 is the starting point
in our proof and is extended via orthogonality arguments.
The basic steps are contained in the following Propositions
2.3-2.5.

In what follows $N$, denotes an integer $\le 10d-1$ and 
$l=[\log_2(\sqrt\la)]=\ell_0$. Then the following estimates hold 
uniformly in $j,k,l$.

\proclaim{Proposition 2.3}

(i) 
For $2^l\le \la^{1/2}$, $k\le j\le l/2$,
$$\aligned
\sup_P\,&\cA_P(2^{j+k-l}, 2^{k-l/2}, 2^{j+k-l}, 2^{k-l/2})
\\&\lc \sup_P\cA_P(2^{-l}, 2^{-j-l/2}, 2^{-l}, 2^{-k-l/2}) +
2^{-l(2d-1)/2}2^{-(j+k)/2}(2^{-2l}\la)^{-N/2}.
\endaligned
\tag 2.17
$$

(ii) For  $k\le j\le \ell_0/2$,
 $$\aligned
\sup_P \,&\cA_P^0(2^{j+k-\ell_0}, 2^{k-\ell_0/2}, 2^{j+k-\ell_0 },
2^{k-\ell_0/2})
\\&\lc
\sup_P\cA_P^0(2^{-\ell_0}, 2^{-j-\ell_0/2}, 2^{-\ell_0}, 2^{-k-\ell_0/2}) +
\la^{-d/2+1/4} 2^{-(j+k)/2}.
\endaligned
\tag 2.18
$$
\endproclaim

\proclaim{Proposition 2.4}

(i)  For $2^l\le \la^{1/2}$, $k\le j\le l/2$,
$$
\aligned\sup_P\,&\cA_P(2^{j-l/2}, 2^{j-l/2}, 2^{k-l/2}, 2^{k-l/2})
\\&\lc
\sup_P\,\cA_P(2^{j+k-l}, 2^{k-l/2}, 2^{j+k-l}, 2^{k-l/2})+
2^{(j+k)(d-1)} 2^k 2^{-l(2d-1)/2}(2^{j+k-2l}\la)^{-N/2}.
\endaligned
\tag 2.19
$$

(ii) For  $k\le j\le \ell_0/2$, $$\aligned\sup_P\,&
\cA_P^0(2^{j-\ell_0/2}, 2^{j-\ell_0/2}, 2^{k-\ell_0/2 }, 2^{k-\ell_0/2})
\\&\lc\sup_P\cA_P^0(2^{j+k-\ell_0}, 2^{k-\ell_0/2}, 2^{j+k-\ell_0},
2^{k-\ell_0/2})+
2^{(j+k)\frac{2(d-1)-N}2 +k}\la^{- \frac{2d-1}4 -\frac N 2 }
...
\endaligned
\tag 2.20
$$
\endproclaim

\proclaim{Proposition 2.5}

(i)  For $2^l\le \la^{1/2}$, $k\le j\le l/2$,
$$
\|T[\si_{j,k,l}]\|\lc
\sup_P\,\cA_P(2^{j-l/2}, 2^{j-l/2}, 2^{k-l/2}, 2^{k-l/2})+
2^{(j+k-l)d/2}(\la 2^{k-3l/2})^{-N/2}.
\tag 2.21$$

(ii) For  $k\le j\le \ell_0/2$, $$\|T[\si^0_{j,k,\ell_0}]\|\lc
\sup_P\,\cA_P^0(2^{j-\ell_0/2}, 2^{j-\ell_0/2}, 2^{k-\ell_0/2 },
2^{k-\ell_0/2})+
2^{(j+k)d/2 - kN/2} \la^{-d/2-N/8}.
\tag 2.22
$$
\endproclaim

Taking these estimates for granted  we can give the
\demo{\bf Proof of Theorem 2.1}

Observe that  since $k\le j\le l/2$ and $2^l\le \la^{1/2}$ the quantities
$2^{-l(2d-1)/2}2^{-(j+k)/2}(2^{-2l}\la)^{-N/2}$,
$2^{(j+k)(d-1)/2 +k-l (2d-1)/2 }(2^{j+k-2l}\la)^{-N/2})$
and
$2^{(j+k-l)d/2}(\la 2^{k-3l/2})^{-N/2}$ are all dominated by a constant
times
$\la^{-(d-1)/2}\min\{
2^{-(l+j+k)/2}, 2^{l/2}\la^{-1/2}\}$, and a combination of  the first parts
of the Propositions 2.3-2.5
gives

$$\align
\big\|T_\la[\si_{j,k,l}]\big\|\lc
& \sup_P\cA_P(2^{-l}, 2^{-j-l/2}, 2^{-l}, 2^{-k-l/2}) +
\la^{-(d-1)/2}\min\{
2^{-(l+j+k)/2}, 2^{l/2}\la^{-1/2}\}.
\endalign
$$
We estimate  the quantities $\cA_P(2^{-l}, 2^{-j-l/2}, 2^{-l}, 2^{-k-l/2})$
  by Proposition 2.2 and (2.5) follows. (2.6) is proved in the same way,
using instead (2.18), (2.20) and (2.22).\qed
\enddemo

\head{\bf 3. Proofs of the Propositions}\endhead

\subhead{Preliminaries}\endsubhead
We begin by stating two  elementary Lemmas which will be
used several times in the proof of Propositions 2.3-5.

\proclaim
{Lemma 3.1} Suppose that
$\zeta\in \cZ_P^{a_P,b_P}(\gamma_1,\gamma_2,\delta_1,\delta_2)$.
Then $\zeta=\sum_{i=1}^M c_i\zeta_i$ where
$\zeta_i\in \cZ_{Q_i} (\eps\gamma_1,\eps\gamma_2,\eps\delta_1,\eps\delta_2)$
with $Q_i\in
B_P^{a_P,b_P}(\gamma_1,\gamma_2,\delta_1,\delta_2)$ and so that
$|c_i|, M\le C_\eps$ (independent of the specific choice of $\zeta$ and
$\gamma$, $\delta$, $P$).
\endproclaim

\demo{Proof} Immediate.\qed \enddemo

\proclaim
{Lemma 3.2}
 Let $Q=(x^Q,z^Q)\in B_P^{a_P,b_P}(\gamma_1,\gamma_2,\delta_1,\delta_2)$
where
$\gamma_1\le \gamma_2\le \eps$,
$\delta_1\le \delta_2\le \eps $.
Suppose that
$0<\widetilde \gamma_1\le \widetilde \gamma_2\le \eps^{-1} \gamma_2$,
$0<\widetilde \delta_1\le \widetilde \delta_2\le \eps^{-1} \delta_2$
and assume that
$$\min\{\frac{\widetilde \gamma_1}{ \widetilde\gamma_2},
\frac{\widetilde \delta_1}{ \widetilde\delta_2}\}\ge \max \{\gamma_2,
\delta_2\}.
\tag  3.1
$$
Then there  are positive constants $C$, $C_1$ (independent of  $\gamma, \widetilde
\gamma,\delta,
\widetilde \delta$, $P$, $Q$) so that for
  $\zeta\in \cZ_Q^{a_P,b_P}(\widetilde \gamma_1,\widetilde
\gamma_2,\widetilde \delta_1,
\widetilde \delta_2)$ the function
$C_1^{-1} \zeta$ belongs to $ \cZ_Q^{a_Q,b_Q}
(C\widetilde \gamma_1,C\widetilde \gamma_2,C\widetilde \delta_1,C\widetilde
\delta_2)$.
\endproclaim

\demo{Proof}
Observe that $$|a_P-a_Q|+|b_P-b_Q|\lc
\max\{\gamma_2,\delta_2\}.$$
The relevant geometry   is then  that by assumption (3.1)  the boxes
$B_Q^{a_P,b_P}(\widetilde \gamma_1,\widetilde \gamma_2,\widetilde \delta_1,
\widetilde \delta_2)$
and
$B_Q^{a_Q,b_Q}(\widetilde \gamma_1,\widetilde \gamma_2,\widetilde \delta_1,
\widetilde \delta_2)$
are contained in fixed dilates of each other. The asserted estimates are
easy to check.\qed
\enddemo

We shall denote by $\eta$ a $C^\infty_0(\bbR)$ function which is supported
in $(-1,1)$ and
satisfies
$\sum_{n\in\bbZ} \eta(\cdot-n)\equiv 1.$
Moreover the $C^\infty_0(\bbR^{d-1})$ function $\chi$ is defined by
$\chi(x_1,\dots,x_{d-1})=\prod_{i=1}^{d-1}\eta(x_i)$.
In the proofs of Propositions 3.3-5 we shall use dilates and translates of
$\eta$ and  $\chi$ to
decompose a suitable cutoff function $\zeta$ as
$$\zeta=\sum_{X,Z\in \bbZ^d}\zeta_{XZ}
\tag3.2$$
the definition of $\zeta_{XZ} $ depends on the particular geometry and is
given by
 (3.12), (3.25) and  (3.30)
below
in the three respective cases.
We shall then employ orthogonality arguments to estimate the operator norm
of
$T_\la[\zeta\si_{j,k,l}]$ in terms of the operator norms of
$$
T_{XZ}:=
T_\la[\zeta_{XZ}\si_{j,k,l}].$$
 This is done by using the Cotlar-Stein Lemma \cite{25, ch. VII,2}.
We  then have to estimate the kernels of
$T^*_{XW} T_{\widetilde X Z}$ and
$T_{XZ} T_{Y\widetilde  Z}^*$.

The kernel of
$T^*_{XW} T_{\widetilde X Z}$  is given by
$$
H(w,z)\equiv H_{XW\widetilde X Z}(w,z)= \int e^{-\ic\la
(\Phi(x,w)-\Phi(x,z))}
\ka_{XW\widetilde XZ}(x,w,z) dx
\tag 3.3
$$
where
$$
\ka_{XW\widetilde XZ}(x,w,z) =\zeta_{\widetilde X Z}(x,z)
\overline{\zeta_{ X W}(x,w)}\si_{j,k,l}(x,z)
\overline{\si_{j,k,l}(x,w)}.
\tag 3.4
$$

The kernel of $T_{XZ}T^*_{Y\widetilde Z}$ is given by
$$
K(x,y)\equiv K_{XZY\widetilde Z}(x,y)
=\int e^{\ic\la(\Phi(x,z)-\Phi(y,z))} \om_{XZY\widetilde Z} (z,x,y) dz
\tag 3.5
$$
with
$$\om_{XZY\widetilde Z}(z,x,y)=\zeta_{ X Z}(x,z)
\overline{\zeta_{ Y\widetilde Z}(y,z)}\si_{j,k,l}(x,z)
\overline{\si_{j,k,l}(y,z)}.
\tag 3.6
$$

Our localizations $\zeta_{XZ}$ will always have the property that
the supports of $\zeta_{XW}$ and $\zeta_{\widetilde X Z}$ are disjoint
whenever
$|X_i-\widetilde X_i|\ge 3$ for some
$i\in \{1,\dots,d\}$. Moreover the supports of
$\zeta_{XZ}$ and $\zeta_{Y\widetilde  Z}$ are disjoint  whenever
$|Z_i-\widetilde Z_i|\ge3$
for some
$i\in \{1,\dots,d\}$.
This implies that
$$
T_{XW}^* T_{\widetilde XZ}=0 \quad \text{ if }|X-\widetilde X|_\infty\ge 3,
\qquad
T_{XZ} T_{Y\widetilde Z}^*=0 \quad \text{ if }|Z-\widetilde Z|_\infty\ge 3.
\tag 3.7
$$

In what follows we shall split variables $X$ and $Z$ in  $\bbZ^{d}$ as
$X=(X', X_d)$,
$Z=(Z', Z_d)$. The geometric meaning of this splitting depends on the
particular situation
 in Propositions 3.3-5.

The main orthogonality properties will always follow from
either the localization properties of the operator in terms of $h$, $V_L h$
or $V_R h$, or by an integration by parts with respect to the directions
orthogonal to $a_P$ or $b_P$.
To describe this we assume
 that $a_P=e_d$, $b_P=e_d$ at a suitable reference point, a situation which
we will
always be able to  achieve by an affine  change of variables as described in
\S2.
If $\Phi_{x'}(x,w)\neq \Phi_{x'}(x,z)$ for all $x$ with $(x,w,z)\in \supp
\ka_{XW\widetilde XZ}$
then we may integrate by parts with respect to the $x'$ variables;
specifically
we have
$$
H(w,z)= (i/\lambda)^N \int
e^{-\ic\la (\Phi(x,w)-\Phi(x,z)}) \cL^N\![\ka_{XW\widetilde XZ}](x,w,z) dx
\tag 3.8
$$
where the differential operator $\cL$ is defined by
$$ \cL g=\text{div}_{x'}\big(\frac{\Phi_{x'}(x,z)-\Phi_{x'}(x,w)}{
|\Phi_{x'}(x,z)-\Phi_{x'}(x,w)|^2}g \big).
\tag 3.9
$$
Similar formulas hold for the $z'$ integration by parts
for the  integral defining $K(x,y)$.

We shall give a proof  of the
 estimates (2.17), (2.19) and (2.21), and the proof of (2.18), (2.20) and
(2.22) is similar.
Here we note that the {\it lower} bound on $|h|$ in the localization
(2.5) is  used in the proof of
estimate (2.14); however it is not needed for the proof  of Propositions
2.3-2.5.

\demo{\bf Remarks on the proof of Proposition 2.2}
In order to prove (2.14) it suffices, by  Lemma 3.1, to  estimate
$\cA_P(\eps 2^{-l}, \eps 2^{-j-l/2}, \eps 2^{-l},\eps 2^{-k-l/2})$
for small $\eps$.

By an affine change of variable as discussed in (2.9)
 we may assume that $P=(0,0)$ and that
$\Phi_{z'x_d}(0,0)=0$, $\Phi_{x'z_d}(0,0)=0$, thus
$\Phi_{x'z'} $ is close to the identity $I_{d-1}$ on the support of $\zeta$
and
the quantities $|\Phi_{z'x_d}(x,z)|$ and $|\Phi_{x'z_d}(x,z)|$ are bounded
by $A_0\eps 2^{-k-l/2}$
for $(x,z)\in \supp \zeta$
(recall that $k\le j$).
 Moreover $a_P=e_d$, $b_P=e_d$;
thus $\zeta$ is, up to a constant, a normalized cutoff function
associated to a box where
$|x'|, |z'|\le \eps 2^{-l}$, $|x_d|\le 2^{-j-l/2}$, $ |z_d|\le 2^{-k-l/2}$.
This puts us in the situation as in the proof of \cite{4, (3.6)}.
If $\cA_P(\eps 2^{-l}, \eps 2^{-j-l/2}, \eps 2^{-l},\eps 2^{-k-l/2})$
 does not vanish identically then the
 function  $|h(x,z)|$ is comparable to $2^{-l}$  on the box
 $B_P^{a_P,b_P}(\eps 2^{-l}, \eps 2^{-j-l/2}, \eps 2^{-l},\eps 2^{-k-l/2})$.
Set $S=T_\la[\zeta\sigma_{j,k,l}]$. For  $j\ge k$  (assumed here)
 the kernel of  $SS^*$ can be estimated using integration by parts, and
all
 the  details of this argument are provided  in \cite{4}.

The estimate (2.15) is more standard, but we sketch the argument for
completeness.  We may assume that $(V_L)_P=\partial_{ z_d}$ and
$(V_R)_P=\partial_{x_d}$ and then ``freezing'' $x_d,z_d$ we may write
$$Sf(x',x_d)=\int_{\Bbb R} S^{x_d,z_d} [f(\cdot,z_d)](x') dz_d.$$
Each $S^{x_d,z_d}$ is an oscillatory integral operator  of the form (1.1)
 in $\Bbb R^{d-1}$ and the mixed Hessian of the
phase function has maximal rank $d-1$;
 however the amplitudes have less favorable differentiability properties.
Note that each $(x',z')$ differentiation
causes a blowup of $O(2^l)=O(\la^{1/2})$.
These estimates for the amplitudes are analogous to the
differentiability properties of symbols of type $(1/2,1/2)$, and in this
situation
the classical bound remains true;  one can
combine   H\"ormander's argument in \cite{16}  with  almost-orthogonality
arguments  in the proof of the
Calder\'on-Vaillancourt theorem for pseudo-differential operators \cite{2}.
See
also  \cite{11} for related but somewhat different arguments for Fourier
integral operators
associated  to canonical graphs. Here it follows that the $L^2$ operator
norm of
$S^{x_d,z_d}$ is $O(\la^{-(d-1)/2})$ uniformly in $x_d,z_d$.
>From the definition of $\sigma_{j,k,l}$ we see that  there are intervals 
$I$ and $J$ of length
$O(2^{-j-l/2})$ and $O(2^{-k-l/2})$, respectively, so that
$S^{x_d,z_d}=0$ unless $x_d\in I$ and $z_d\in J$. Thus from
applications of Minkowski's and
Cauchy-Schwarz' inequalities it follows that $\|S\|\lc 2^{-(j+l/2)/2}
2^{-(k+l/2)/2}
\la^{-(d-1)/2}$.
(2.16) is proved in the same way.\qed

\enddemo

\demo{ \bf Proof of Proposition 2.3} Fix $P$.
By Lemma 3.1
it suffices  to estimate
$\|T_\la[\zeta \si_{j,k,l}]\|$ where $\zeta$ belongs to $\cZ_P^{a_P,b_P}(\eps
2^{j+k-l},
\eps 2^{k-l/2},\eps 2^{j+k-l},\eps 2^{k-l/2})$, with norm independent of
$P$.

By an affine change of variable as discussed in (2.9)
 we may assume that $P=(0,0)$ and that
$\Phi_{z'x_d}(0,0)=0$, $\Phi_{x'z_d}(0,0)=0$, hence
$$
\gather
\|\Phi_{x'z'}-I\|\le 2^{-d}
\tag 3.10\\
|\Phi_{z'x_d}(x,z)|+|\Phi_{x'z_d}(x,z)|\le A_0\eps 2^{k-l/2}
\tag 3.11
\endgather
$$
for $(x,z)\in \supp \zeta$. Moreover $a_P=e_d$, $b_P=e_d$;
thus $\zeta$ is, up to a constant, a normalized cutoff function
associated to a box where
$|x'|, |z'|\le \eps 2^{j+k-l}$, $|x_d|, |z_d|\le 2^{k-l/2}$.

For $X, Z\in\bbZ^d$ let
$$\multline
\zeta_{XZ}(x,z)=\\\zeta(x,z)
\chi(2^l\eps^{-1} x'-X')\eta( 2^{j+l/2}\eps^{-1}x_d-X_d)
\chi(2^l\eps^{-1} z'-Z')\eta( 2^{k+l/2}\eps^{-1}z_d-Z_d)
\endmultline
\tag 3.12
$$
and let $T_{XZ}=T_\la[\zeta\sigma_{j,k,l}] $.
By (3.10-11)
and
Lemma 3.2 there are positive constants $C, C_1$ so that
$C_1^{-1} \zeta_{XZ}$ belongs to
$\cZ_Q^{a_Q,b_Q}(C2^{-l},C 2^{-j-l/2},C 2^{-l},C2^{-k-l/2})$.
Thus
$$\|T_{XZ}\|\lc \sup_P\,\cA_P(2^{-l}, 2^{-j-l/2}, 2^{-l},2^{-k-l/2})
\tag 3.13
$$
and it remains to show almost orthogonality of the pieces $T_{XZ}$.

By our localization the orthogonality properties  (3.7) are satisfied.
Therefore the assertion (2.17) follows from
$$\aligned
\big\|T_{XW}^* &T_{\widetilde XZ}\big\|
\lc 2^{-l(2d-1)-j-k} (\la2^{-2l}|W'-Z'|)^{-N} \\
&\text{ if }
2A_0 |W'-Z'|\ge |W_d-Z_d| \text{ and } |W'-Z'|\ge C_1,
\endaligned
\tag 3.14$$
$$
T_{XW}^* T_{\widetilde XZ}=0,\quad \text{ if }
2A_0|W'-Z'|< |W_d-Z_d|  \text{ and } |W_d-Z_d|\ge C_1,
\tag 3.15
$$
(for suitable $C_1\gg1 $)
and
$$\aligned
\big\|T_{XZ} &T_{Y\widetilde Z}^*\big\|
\lc 2^{-l(2d-1)-j-k} (\la2^{-2l}|X'-Y'|)^{-N}, \\
&\text{ if }
2A_0 |X'-Y'|\ge |X_d-Y_d| \text{ and } |X'-Y'|\ge C_1,
\endaligned
\tag 3.16
$$
$$
T_{XZ} T_{Y\widetilde Z}^*=0, \quad\text{ if }
2A_0|X'-Y'|< |X_d-Y_d|\text{ and } |X_d-Y_d|\ge C_1.
\tag 3.17
$$

We now show (3.15) and (3.14).
The kernel $H$ of
$T_{XW}^* T_{\widetilde XZ}$
 is given by (3.3), (3.4).
 In order to see (3.15) pick  points
$(x,w)\in \supp\zeta_{XW}$ and
$(x,z)\in \supp\zeta_{\widetilde XZ}$
and also assume that $(x,w)$ and $(x,z)$ belong to $\supp \si_{j,k,l}$ (if
there
are no two such points then
$T_{XW}^* T_{\widetilde XZ}=0$).
By definition of $\si_{j,k,l}$
we have
$$|h(x,z)-h(x,w)|\le 2^{-l+2}.
\tag 3.18
$$
Also for all $ (x,\tilde z)\in \supp \zeta$ we have that
$$|h_{z_d}(x,\tilde z)-V_L h(x,\tilde z) |\le A_0\eps 2^{k-l/2}
\tag 3.19$$
so that $|h_{z_d}(x,\tilde z)|\ge 2^{k-l/2-2}$.
Note that $\eps 2^{-k-l/2}(|W_d-Z_d|-2)
\le|w_d-z_d|\le
2^{-k-l/2}(|W_d-Z_d|+2)\eps$
and
 $|w'-z'|\le C_0(|W'-Z'|+2) 2^{-l}\eps$. Therefore
$$
\align
|h(x,w)-h(x,z)|&\ge |h_{z_d}(x,z)||w_d-z_d|-
A_0|w'-z'|
\\
&\ge 2^{k-l/2-4} \eps 2^{-k-l/2}|W_d-Z_d|=\eps 2^{-l-2}|W_d-Z_d|
\tag 3.20
\endalign
$$
if
$|W_d-Z_d|\ge \max\{2A_0 |W'-Z'|_\infty, 2A_0\}$ and  $|W_d-Z_d|\ge 2A_0$.
Observe that (3.18) and (3.20) can hold simultaneously only when
$|W_d-Z_d| $ stays bounded; this implies (3.15).

Now assume that
$|W_d-Z_d|\le 2A_0 |W'-Z'|_\infty$, and we show (3.14) if $|W'-Z'|\ge C_1$
for sufficiently large $C_1$.

We perform
 integration by parts with respect to the $x'$ variables in (3.3),
using  (3.8/3.9).
Now in view of
(3.10/11) we have
$$\align
|\Phi_{x'}(x,w)-\Phi_{x'}(x,z)|
&\ge |\Phi_{x'z'}(x,z)(w'-z')|- A_0\eps 2^{k-l/2}|w_d-z_d|- A_0|w-z|^2
\\
&\ge \eps 2^{-l}|W'-Z'|
\tag 3.21
\endalign
$$
if $|W'-Z'|\ge C_1$ for suitable $C_1$. Moreover
 the $x$-derivatives of
$\Phi(x,w)-\Phi(x,z)$ are  $O(2^{-l}|W'-Z'|)$, and differentiating the
symbol causes a blowup of $O(2^{l})$ for each differentiation.
Thus for $|W'-Z'|\ge C_1$
$$|\cL^N(\kappa_{XW\widetilde X Z})|\lc
(\la 2^{-2l}|W'-Z'|)^{-N}.$$
Taking into account the $x$ support this yields the estimate
$$|H(w,z)|\lc 2^{-l(d-1)}2^{-l/2-j}(\la 2^{-2l}|W'-Z'|)^{-N}.
$$
By Schur's test we have to bound
$\sup_w\int|H(w,z)|dz$ and $\sup_z\int|H(w,z)|dw$. Since the integrals
are extended over sets of measure
$O(2^{-l(d-1)}2^{-l/2-k})$ we obtain the bound (3.14).

We still have to estimate the kernel $K$ given by (3.5), (3.6).
Note that
$|h_{x_d}(x,\tilde z)-V_R h(x,\tilde z) |\le A_0\eps 2^{k-l/2}$ so that
$|h_{x_d}(x,z)|\approx 2^{j-l/2}$ (recall that $j\ge k$). Thus in place of
(3.20) we have
$$|h(x,z)-h(y,z)|\ge 2^{j-l/2}|x_d-y_d|-
A_0|x'-y'| \tag 3.22
$$
and in place of (3.21) we have
$$
|\Phi_{z'}(x,z)-\Phi_{z'}(y,z)|
\ge |\Phi_{z'x'}(y,z)(x'-y')|-A_0\eps 2^{k-l/2}|x_d-y_d|- A_0|x-y|^2.
\tag 3.23
$$
Since $|x_d-y_d|\approx |X_d-Y_d|$ we proceed as before to obtain  (3.16)
and (3.17).\qed
\enddemo

\demo{ \bf Proof of Proposition 2.4}
We continue to use the same notations as in the previous proof
although our localizations are with respect to different
(larger) boxes.
By Lemma 3.1
it suffices  to estimate  the operator norm of
$T_\la[\zeta \si_{j,k,l}]$ where now $\zeta\in \cZ_P(\eps 2^{j-l/2},
\eps 2^{j-l/2},\eps 2^{k-l/2},\eps 2^{k-l/2})$. Again we may assume that by
an affine
 change of variable $P=(0,0)$ and that
$\Phi_{z'x_d}$, $\Phi_{x'z_d}$ vanish at $(0,0)$.
It follows that
$$
|\Phi_{z'x_d}(x,z)|+|\Phi_{x'z_d}(x,z)|\le A_0\eps 2^{j-l/2}, \quad
(x,z)\in \supp \zeta,
\tag 3.24
$$
and again $a_P=e_d$, $b_P=e_d$.

For $X, Z\in\bbZ^d$ we now define
$$
\multline
\zeta_{XZ}(x,z)= \zeta(x,z) \times\\
\chi(2^{-j-k+l}\eps^{-1} x'-X')\eta( 2^{-k+l/2}\eps^{-1}x_d-X_d)
\chi(2^{-j-k+l}\eps^{-1} z'-Z')\eta( 2^{-k+l/2}\eps^{-1}z_d-Z_d)
\endmultline
\tag 3.25
$$
and set $T_{XZ}=
T_\la[\zeta_{XZ} \si_{j,k,l}]$. In view of  (3.24), Lemma 3.1 and Lemma 3.2
$$\|T_{XZ}\|\le \sup_P\,\cA_P(2^{j+k-l}, 2^{k-l/2}, 2^{j+k-l},2^{k-l/2}).$$
To show the orthogonality observe that (3.7) remains valid.
Moreover the width of the smaller boxes in the $z_d$ direction is comparable
to the $z_d$-width
of the original  boxes, namely $\approx 2^{k-l/2}$.
This shows that
$$
T_{XW}^* T_{\widetilde X Z}=0\qquad\text{ if } |W_d-Z_d|\ge C_1
\tag 3.26
$$ for sufficiently large $C_1$.

This estimate is complemented by
$$
\big\|T_{XW}^* T_{\widetilde XZ}\big\|
\lc 2^{2(j+k)(d-1)}2^{2k}2^{-l(2d-1)} (\la 2^{j+k-2l}|W'-Z'|)^{-N},
\quad\text{ if }
|W'-Z'|\ge C_1,
\tag 3.27$$
for large $C_1$.

To see (3.27) we integrate by parts with respect to $x'$.
Our kernel is still given by (3.3), (3.4).
To perform the integration by parts we may assume that
$|W_d-Z_d|\le C_1$ by (3.26). We now see from (3.24) that
$$
|\Phi_{x'}(x,w)-\Phi_{x'}(x,z)|
\ge |\Phi_{x'z'}(x,z)(w'-z')|-A_0\eps 2^{j-l/2}|w_d-z_d|- A_0|w-z|^2
$$
but $|w'-z'|\approx  |W'-Z'| \eps2^{j+k-l}$, and
$|w_d-z_d|\le 2^{k-l/2}\eps|W_d-Z_d|\le C\eps 2^{k-l/2}$. Thus if $|W'-Z'|$
is sufficiently large
we have the lower bound
$$
|\Phi_{x'}(x,w)-\Phi_{x'}(x,z)|
\gc 2^{j+k-l}|W'-Z'|
$$ for $(x,w,z)\in \supp \ka_{XW\widetilde XZ}$.
Therefore analyzing
$\cL^N(\kappa_{XW\widetilde X Z})$ as in the proof of Proposition 3.3 we
see that
$$|\cL^N(\kappa_{XW\widetilde X Z})|\lc ( 2^{j+k-2l} |W'-Z'|)^{-N}.$$ From this
 we get the pointwise bound
$$|H(w,z)|\lc 2^{(j+k-2l)(d-1)}2^{k-l/2}
(\la 2^{j+k-2l}|W'-Z'|)^{-N}.$$
For Schur's test we have to integrate this in $x$ or $y$  over a set of
measure
$2^{(j+k-l)(d-1)}2^{k-l/2}$ and we obtain in fact a slightly better estimate than (3.27).

Next,  it remains to show that
$$
\aligned
\big\|T_{XZ} &T_{Y\widetilde Z}^*\big\|
\lc 2^{2(j+k)(d-1)}2^{2k}2^{-l(2d-1)} (\la 2^{j+k-2l}|X'-Y'|)^{-N} ,
\\&\text{ if }
2A_0 |X'-Y'|\ge |X_d-Y_d| \text{ and } |X'-Y'|\ge C_1,
\endaligned
\tag 3.28
$$
and
$$
T_{XZ} T_{Y\widetilde Z}^*=0 \quad\text{ if }
2A_0|X'-Y'|< |X_d-Y_d|\text{ and } |X_d-Y_d|\ge C_1.
\tag 3.29
$$

The proof of these  estimates is  similar to the proof of
the corresponding estimates in
Proposition 2.3.
The estimate (3.22) continues to hold and the estimate (3.23) is replaced by
the weaker estimate
$$|\Phi_{z'}(x,z)-\Phi_{z'}(y,z)|
\ge |\Phi_{z'x'}(y,z)(x'-y')|-A_0\eps 2^{j-l/2}|x_d-y_d|- A_0|x-y|^2
$$
which however still gives the asserted bound since
$|x_d-y_d|\approx 2^{k-l/2}\eps |X_d-Y_d|$ and
$|x'-y'|\approx 2^{-l}|X'-Y'|$. \qed
\enddemo

\demo{\bf Proof of Proposition 2.5}
We may assume that the support of $\zeta$ is small ({\it i.e.} contained in
a ball of radius $\eps$).
By Lemma 3.1
it suffices  to estimate  the operator norm of
$T_\la[\zeta \si_{j,k,l}]$ where now $\zeta\in \cZ_P(\eps,
\eps,\eps,\eps)$. By  affine changes of variables we may assume that
$P=(0,0)$ and that
$\Phi_{z'x_d}$, $\Phi_{x'z_d}$ vanish at $(0,0)$.  Thus
$$|\Phi_{z'x_d}(x,z)|+|\Phi_{x'z_d}(x,z)|\le A_0\eps, \qquad
(x,z)\in \supp \zeta.
$$

For $X, Z\in\bbZ^d$ we now consider $T_{XZ}=T_\la[\zeta\si_{j,k,l}]$ with
$$
\multline
\zeta_{XZ}(x,z)=\zeta(x,z)
\chi(2^{-j+l/2}\eps^{-1} x'-X')
\times \\
 \eta( 2^{-j+l/2}\eps^{-1}x_d-X_d)
\chi(2^{-k+l/2}\eps^{-1} z'-Z')\eta( 2^{-k+l/2}\eps^{-1}z_d-Z_d),
\endmultline
\tag 3.30
$$
and again $\|T_{XZ}\|\lc\sup_P\,
\cA_P(2^{j+k-l}, 2^{k-l/2}, 2^{j+k-l},2^{k-l/2}).$

%

For the orthogonality of the pieces  we now use besides (3.7)
the assumptions (2.3). By our choice of $\eps$
we have that
$$|V_L^2 h-\partial_{z_d} V_L h|\le A_0\eps \le A_1/10$$
and similarly
$$|V_R^2 h-\partial_{x_d} V_R h|\le  A_1/10.$$
Thus
$$
\align
4\cdot 2^{k-l/2}&\ge |V_L h(x,w)-V_L h(x,z)|\ge (2A_1)^{-1}
|w_d-z_d|-A_0|w'-z'|
\\
4\cdot
 2^{j-l/2}&
\ge |V_R h(x,z)-V_R h(y,z)|\ge (2A_1)^{-1} |x_d-y_d|-A_0|x'-y'|.
\endalign
$$
This shows that
$$
\align
 T_{XW}^* T_{\tilde XZ}&=0 \text{ if }|W_d-Z_d|\ge C
\tag 3.31
\\
T_{XZ} T_{Y\tilde Z}^*&=0 \text{ if }|X_d-Y_d|\ge C
\tag 3.32
\endalign
$$
for $C=10 A_0A_1$.

Now assume that
$|W_d-Z_d|\le C_1$. Then
if $(x,w,z)\in \supp \overline{\zeta_{XW}}\zeta_{XZ}$ we have
$$
|\Phi_{x'}(x,w)-\Phi_{x'}(x,z)|
\ge |\Phi_{x'z'}(x,z)(w'-z')|-A_0\eps |w_d-z_d|- A_0|w-z|^2
$$
but now  $|w'-z'|\approx  |W'-Z'| 2^{k-l/2}$, and
$|w_d-z_d|\le 2 C 2^{k-l/2}$. Thus for large  $|W'-Z'|$
we have the lower bound
$$
|\Phi_{x'}(x,w)-\Phi_{x'}(x,z)|
\ge 2^{k-l/2}|W'-Z'|
$$ and it follows that
$$|\cL^N[\kappa_{XW\widetilde X Z}](x,w,z)|\lc
(\la 2^{k-3l/2}|W'-Z'|)^{-N}.$$
Consequently $|H(w,z)|\lc 2^{(j-l/2)d}
(\la 2^{k-3l/2}|W'-Z'|)^{-N}.$
To apply Schur's test we observe that for fixed
$z$ the $w$ integral is extended over a set of measure
$O(2^{(k-l/2)d})$ (likewise for fixed $w$ the $z$ integral). We obtain
the bound
$$
\|T_{XW}^* T_{\tilde XZ}\|\lc
2^{(j+k-l)d} (\la 2^{k-3l/2}|W'-Z'|)^{-N}
\tag 3.33
$$
if $|W'-Z'|\ge C'$.
By a similar argument
$$
 \|T_{XZ} T_{Z\tilde Z}^*\|\lc
2^{(j+k-l)d} (\la 2^{j-3l/2}|X'-Y'|)^{-N}.
\tag 3.34
$$
The asserted estimate (2.21) now follows from combining
 (3.31-34) and the estimate for the individual pieces.\qed
\enddemo

\head{ \bf 4. One-sided type three singularities} \endhead
In this section we discuss the  proof of Theorem 1.2. The reasoning
 is very close to the one given by the authors in \cite{9}, but the
assumptions there are somewhat different. We  thus only
sketch the proof and refer the reader to \cite{8}, \cite{9} for details of
some of the  arguments.

First we shall  need  an  extension of
Theorem 1.1 to oscillatory integral operators of the form
$$\cT_\mu f(x)= \int f(y) \int e^{\ic\mu\psi(x,y,\vth)} a(x,y,\vth) d\vth
\,dy
$$
where the frequency variable $\vth$ lives  in an open set $\Theta\subset
\bbR^N$ and we assume that $a\in C^\infty_0(\Om_L\times\Om_R\times \Theta)$.
It is assumed that  $\psi$ is a
nondegenerate phase function in the sense of H\"ormander \cite{15} (but
not necessarily homogeneous),
{\it i.e.} $\nabla_{\vartheta_i}(\nabla_{x,y,\vth}\Psi)$, $i=1,\dots, N$ are
linearly independent.
The canonical relation $C_\psi\subset T^*\Om_L\times T^*\Om_R$ is given by
$$C_\psi=\{(x,\psi_x,y,-\psi_y):\psi_\vth=0\}.$$

\proclaim{Lemma 4.1} Suppose that the projections
 $\pi_L:C_\psi\to T^*\Om_L$,
 $\pi_R:C_\psi\to T^*\Om_R$ are of type $\le 2$. Then $\|\cT_\mu\|_{L^2\to
L^2}
=O(\mu^{-(d+N-1)/2-1/4})$, $\mu\to\infty$.
This  estimate
is stable under small perturbations of $\psi$ and $a$ in the
$C^{\infty}$-topology.
\endproclaim

The reduction to the situation in Theorem 1.2 involves  canonical
transformations on $T^*\Om_L$ and
 $T^*\Om_R$ and then as in \cite{15} an application of the method of
 stationary phase
to reduce the number of frequency variables (see \cite{8} for details).

The following Lemma deals with phase functions $\Phi(x,z)$ without
frequency variables.

\proclaim{Lemma 4.2} Let $\Phi$ be a real--valued phase function defined
near
$(x^0,z^0)$ and assume that
$\nabla_x (\det \Phi_{xz}(x^0,z^0))\neq 0$, and
$|\det \Phi_{xz}(x^0,z^0)|+|V_R \det \Phi_{xz}(x^0,z^0)|\le c
|\nabla_x (\det \Phi_{xz}(x^0,z^0))|$. Let
$M>0$.

Then, if $c$ is sufficiently small,
 there  are  neighborhoods $\Omega_L^0$ of $x_0$, $\Omega_R^0$ of $z_0$,
neighborhoods $\cU$ and $\cV$ of $(x_0,\nabla_x\Phi(x^0,z^0))$ in
$T^*\Omega_L$,
a canonical transformation $\chi:\cU\to \cV$, and a unitary operator
$U_\la$, so that the following statements hold if $\sigma$ is supported in
$\Om_L^0\times \Om_R^0$.

(i) If $T_\la$ is the integral operator with kernel
$\sigma(x,z) e^{\ic\la\Phi(x,z)}$ then
$$U_\la T_\la=S_\la+R_\la$$
where
 $S_\la$ is an integral operator with kernel
$\tau(x,z) e^{\ic\la\Psi(x,z)}$ and $\|R_\la\|_{L^2\to L^2}=O(\la^{-M})$,

(ii) If $\cC_\Phi=\{(x,\Phi_x,z,-\Phi_z), (x,z)\in \supp \sigma\}$
then for
$\cC_\Psi=\{(x,\Psi_x,z,-\Psi_z), (x,z)\in \supp \tau\}$
we have
$$
\cC_\Psi\subset \{(\chi(x,\xi), z, \zeta): (x,\xi,z,\zeta)\in \cC_\Phi\}.
$$

(iii) $\nabla_z (\det \Psi_{xz})\neq 0$ for $(x,z)\in \supp \tau$.
\endproclaim
\demo{Proof}
This can be extracted from the arguments in \S4 of \cite{9}.
\enddemo
\demo{\bf  Proof of Theorem 1.2}
We work with $T_\la$ as in (1.1) where  $(x,z)$ is close to the origin, and
the origin lies on the singular surface $\{(x,z):\det \Phi_{xz}=0\}$.
We may assume, after a change of variable in $z$ that
$$
\Phi(0,z)=0,\qquad \Phi_{x'z'}(0,z)=I, 
\qquad
\Phi_{x'z_d}(0,z)=0
\tag 4.1
$$
(\cf. the proof of  Lemma 2.7 in \cite{9})
and by a change of variable in $x$ we may also assume
$$\Phi_{z'x_d}(0,0)=0.$$
We assume that $\pi_L$ is of type $\le 2$ and that $\nabla_{x,z}(\det
\Phi_{xz})\neq 0$ where $\det \Phi_{xz}$ vanishes.

If  $\Phi_{x_dx_dz_d}\neq 0$ or $\Phi_{x_dz_dz_d}\neq 0$
 or
$\Phi_{x_dz_dz_dz_d}\neq 0$
then we have a fold or cusp singularity and
 better results then the one claimed in Theorem 1.2 were proved in \cite{8},
\cite{9}.
Therefore assume that
$\Phi_{x_dx_dz_d}$, $\Phi_{x_dz_dz_d}$ and
$\Phi_{x_dz_dz_dz_d}$ are small.
By (4.1) and  Lemma 4.2 we may
assume the more restrictive assumption that
$\nabla_z(\det\Phi_{xz})\neq 0$ which near the origin is equivalent
with $\nabla_z(\Phi_{x_dz_d})\neq 0$, again by (4.1).
After  a rotation we may assume
$$
\Phi_{x_dz_d z_1}(0,0)\neq 0.
\tag 4.2
$$

We now consider the  operator $T_\la T_\la^*$ and estimate it by the
slicing technique in \cite{9} (also familiar from the proof of
Strichartz estimates). Now $T_\la T_\la^* f(x)= \int
\cK^{x_d,y_d}[f(\cdot,y_d)](x') dy_d$ where the kernel of
$\cK^{x_d,y_d}$ as an integral operator acting on functions in $\bbR^{d-1}$
is given by
$$
K^{x_d,y_d}(x',y')=\int e^{\ic\la[\Phi(x',x_d,z)-\Phi(y',y_d,z)]}
\sigma(x,z)
\overline{ \sigma(y,z)} dz.
$$
The computation in \cite{9} shows that after  rescaling the estimation is
reduced to showing that two   integral operators $\cH_\mu^\pm
\equiv \cH_{\mu,\gamma,c}^\pm$ with kernel
$$H_\mu^\pm(u,v)=\int e^{\ic \Psi^\pm(u,v,z;\gamma,c)} b_{\gamma,c}(u,v,z)
dz
$$
are bounded on $L^2(\Bbb R^{d-1})$ with norm
$O(\mu^{-(d-1)-1/4})$, $\mu\ge 1$. Here $b_{\gamma,c}$ is $C^\infty_0$
and $$\Psi^\pm(u,v,z;\gamma,c)=\inn{u-v}{\Phi_{x'}(0,c,z)}\pm
\Phi_{x_d}'(0,c,z)+
r_{\pm}(u,v,z,\gamma,c)
$$
where $\gamma$ and $c$ are small  parameters and
$r_{\pm}(u,v,z,0,c)\equiv 0$. The dependence of
 $\Psi^\pm$, $r_\pm$ and $b$ on $\gamma$ and $c$
 is smooth and the bounds have to be uniform for small $\gamma,c$. For this
it
 remains to show that the operators $\cH_\mu^{\pm}$ are oscillatory
integral operators with two-sided type two singularities to which we can
apply Lemma 4.1 in
 $d-1$ dimensions (with $d$ frequency variables). It suffices to check the
type two condition
at $\gamma=0$, $c=0$.

The condition (4.2) guarantees that $\Psi^\pm$  is indeed a nondegenerate
phase function
with critical set
$$\align\text{Crit}_{\Psi^\pm}&=\{(u,v,z):\nabla_z\Psi^\pm=0\}\\
&=\{(u,v,z):v=u\mp\Phi^{x'z'}\Phi_{z'x_d},
\quad\Phi_{x_dz_d}-\Phi_{x_dz'}\Phi^{z'x'}\Phi_{x'z_d}=0\}
\endalign
$$
where the $\Phi$ derivatives are evaluated at $(0,z)$.
At $x=0$ the second equation becomes
$\Phi_{x_dz_d}(0,z)=0$  and by (4.2) we may solve this equation expressing
$z_1$
as a function $z_1^\pm$ of $\tilde z=(z'',z_d)$ with
$z''=(z_2,\dots,z_{d-1})$.
Set $G^\pm(\tilde z)=\Phi_{x'}(0,z_1^{\pm}(\tilde z), \tilde z)$ and
$B^\pm(\tilde z)=
\Phi^{x'z'} (0,z_1^\pm(\tilde z), \tilde z)
\Phi_{z'x_d}(0,z_1^\pm(\tilde z), \tilde z)
$
then the canonical relation for vanishing  $\gamma, c$ is given by
$$C_{\Psi^\pm}=\{(u, G^\pm(\tilde z), u\mp B^\pm(\tilde z), G^\pm(\tilde
z))\}
$$
which is parametrized by the coordinates $(u,\tilde z)$. The derivative 
 of
the projection to $T^*\Om_L$  in these coordinates is given by
$$
 DG^\pm= 
\pmatrix \Phi_{x'z_1} \frac{\partial z_1}{\partial \tilde z}+\Phi_{x'\tilde z}\endpmatrix
$$ and by (4.1) we see that its determinant equals
$(-1)^{d-1} \partial z_1^\pm/\partial z_d$ and 
$\widetilde V_L=\partial/\partial_{z_d}$ is a kernel vector field for the left
projection.
Moreover
 $\widetilde V_R=\widetilde V_L+\sum_{i=1}^{d-1} c_i(z)\partial/\partial_{u_i}$ so
that
$\widetilde V_L$ and $\widetilde V_R$ coincide when acting on the
determinant.
By implicit differentiation we see that
$\partial_{z_d}^kz_1^\pm - \Phi_{z_dx_dz_1}^{-1}\Phi_{x_dz_d^{k+1}}$
 belongs to the ideal generated by $\Phi_{x_dz_d^{j}}$, $j\le k$. 
>From this one
deduces that $\widetilde V_L$ and $\widetilde V_R$ are of type
$\le k-1$ if one of the derivatives  $\Phi_{x_dz_d^j}$ , $j\le k+1$,
does not vanish. We apply this for $k=3$
to conclude the proof.\qed

\enddemo

\enddocument